\documentclass[sn-mathphys-num]{sn-jnl}

\usepackage{graphicx}%
\usepackage{multirow}%
\usepackage{amsmath,amssymb,amsfonts}%
\usepackage{amsthm}%
\usepackage{mathrsfs}%
\usepackage[title]{appendix}%
\usepackage{xcolor}%
\usepackage{textcomp}%
\usepackage{manyfoot}%
\usepackage{booktabs}%
\usepackage{algorithm}%
\usepackage{algorithmicx}%
\usepackage{algpseudocode}%
\usepackage{listings}%

\usepackage{cjhebrew}

\usepackage{float}

\usepackage{subfigure}

\usepackage{polynom}

\usepackage[T1]{fontenc}

\usepackage{ulem}

\usepackage{natbib}

\usepackage{imakeidx}

\usepackage{fancybox}

\usepackage{fancyhdr}

\usepackage[nice]{nicefrac}

\usepackage{url}

\usepackage{pgfplots}
\pgfplotsset{trig format plots=rad}
\pgfplotsset{compat=1.8}

\usepackage{tikz}
\usetikzlibrary{positioning,chains,fit,shapes,calc,arrows,patterns}
\usetikzlibrary{calc,patterns,angles,quotes}
\usetikzlibrary{decorations.pathreplacing}
\usepgfplotslibrary{fillbetween}
\usepackage{filecontents}
\usetikzlibrary{pgfplots.dateplot}
\usepackage{pgfplotstable}

\theoremstyle{thmstyleone}%
\newtheorem{theorem}{Theorem}
\newtheorem{proposition}[theorem]{Proposition}%

\theoremstyle{thmstyletwo}%
\newtheorem{example}{Example}%
\newtheorem{remark}{Remark}%

\theoremstyle{thmstylethree}%
\newtheorem{definition}{Definition}%

\begin{document}

\title[Article Title]{Mathematics of Family Planning in Talmud}


\author[1]{\fnm{Simon} \sur{Blatt}}\email{simon.blatt@plus.ac.at}

\author[2]{\fnm{Uta} \sur{Freiberg}}\email{uta.freiberg@mathematik.tu-chemnitz.de}

\author*[2]{\fnm{Vladimir} \sur{Shikhman}}\email{vladimir.shikhman@mathematik.tu-chemnitz.de}

\affil[1]{\orgname{Paris Lodron Universit\"at Salzburg}, \orgaddress{\street{Hellbrunner Str. 34}, \city{Salzburg}, \postcode{5020}, \country{Austria}}}

\affil[2]{\orgname{Chemnitz University of Technology}, \orgaddress{\street{Reichenhainer Str. 41}, \city{Chemnitz}, \postcode{09126}, \country{Germany}}}


\abstract{Motivated by the commitments from the Talmud in Judaism, we consider the family planning rules which require a couple to get children till certain numbers of boys and girls are reached. For example, the rabbinical school of Beit Hillel says that one boy and one girl are necessary, whereas Beit Shammai urges for two boys. Surprisingly enough, although the corresponding average family sizes differ in both cases, the gender ratios remain constant. We show more that for any family planning rule the gender ratio is equal to the birth odds. The proof of this result is given by using different mathematical techniques, such as induction principle, Doob's optional-stopping theorem, and brute-force. We conclude that, despite possible asymmetries in the religiously motivated family planning rules, they discriminate neither boys nor girls. }




\maketitle

\section{Introduction}
\label{intro}
The Talmud is one of the most important scriptural sources in Judaism and lays the basis for the Jewish religious laws. The noun ''Talmud'' is derived from the Hebrew verb root \<lmd>, meaning ''to learn''. The latter denotes both the activity of studying and teaching, as well as the object of study. Talmud has two components: the Mishnah, completed around 200 AD, a written compendium of the oral Torah with the canonical collection
of Jewish laws, and the Gemara, with the discussions about these laws, which were recorded in Babylonia around 500 AD.
The teaching of material in the Talmud is presented in the form of a dialog between different rabbinical doctrines in order to come to a decision at the end. This reflects the prevailing state of the tradition and shows how the Torah rules were understood and interpreted by the rabbis in practice and in everyday life.
To make the presentation accessible, the Romm brothers designed the meanwhile classical Talmud leaflet in their publishing house in Vilnius, Lithuania, 1880-1886. Typically, in the middle of a page there is the Mishnah part to be discussed. From line 14 the corresponding Gemara part is written in square letters. To the right one sees the commentary of maybe the ultimately most influential rabbi  Rashi, which stands for Rabbi Schlomo ben
Jizchak, who lived and worked in France and Germany in the 11th century. On the left there are commentaries given by his students. All around there are smaller commentaries attributed to other distinguished scholars. 

The overall focus on discussion in the Jewish religious thoughts historically gave rise to rabbinical schools that promote their way of interpretation for centuries. Two main rabbinical schools, whose opinions on the family planning happen to be relevant also for us here, are attributed to Hillel and Shammai. Hillel was one of the most important Pharisee rabbis from the time before the destruction of the Second Temple in 70 AD. Being the head of the Sanhedrin, the highest Jewish religious and political authority of that time and also the supreme court, he also founded a school of biblical exegesis, known now as the House of Hillel or Beit Hillel in Hebrew transliteration. Shammai was a Jewish scholar of the 1st century and the most important contemporary of Hillel, although he was around 60 years younger. In contrast to the mild Hillel, Shammai was considered strict and irritable in religious and ethical matters. In most cases, though not always, House of Hillel's opinion is the more lenient and tolerant than the House of Shammai's one. While the terms ''liberal'' and ''conservative'' may not perfectly capture the nuances of their positions, Hillel is generally considered to have been more flexible in his interpretations of Jewish law compared to Shammai. On the other hand, Shammai tended to be more stringent in his interpretations, prioritizing strict adherence to the law. Let us cite on a typical controversy between the both Houses, see Talmud, Shabbat 31a:6:

\begin{quote}
    {\it There was another incident involving one gentile who came before Shammai and said to Shammai: Convert me on condition that you teach me the entire Torah while I am standing on one foot. Shammai pushed him away with the builder’s cubit in his hand. This was a common measuring stick and Shammai was a builder by trade. The same gentile came before Hillel. He converted him and said to him: That which is hateful to you do not do to another; that is the entire Torah, and the rest is its interpretation. Go study.}
\end{quote}
In Jewish tradition, however, the disagreement of Hillel and Shammai is viewed as one which had lasting positive value, see Talmud, Pirkei Avot 5:17: 

\begin{quote}
    {\it Every dispute that is for the sake of Heaven, will in the end endure; But one that is not for the sake of Heaven, will not endure. Which is the controversy that is for the sake of Heaven? Such was the controversy of Hillel and Shammai. And which is the controversy that is not for the sake of Heaven? Such was the controversy of Korah and all his congregation.}
\end{quote}

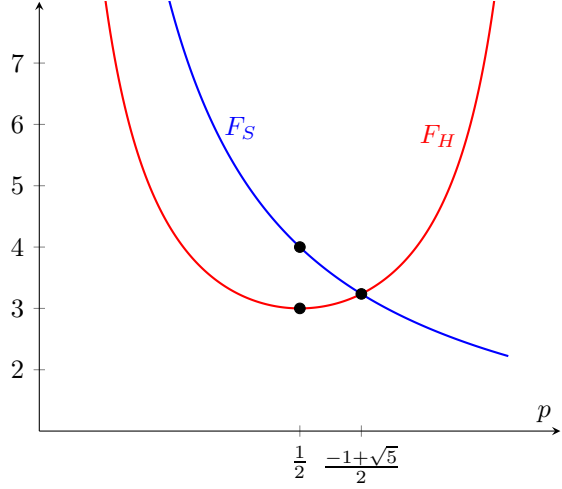
\begin{figure}[ht]
\centering
\begin{tikzpicture}[style=thick]
\begin{axis}
[
  x axis line style = {-latex},
  y axis line style = {-latex},
  axis lines = center,
  xmin=0, xmax=1, ymin=1, ymax=8,
  xlabel = {$p$},
  ylabel = {},
  xtick={0.5,0.618034}, 
  xticklabels={$\frac{1}{2}$,$\frac{-1+\sqrt{5}}{2}$},
  yticklabels={2,3,4,5,6,7},
  ytick={2,3,4,5,6,7}
]
\addplot[thick,
    domain=0.1:0.9, 
    samples=200, 
    color=red,
]
{(x^2 - x + 1)/(x-x^2)} node[left,pos=0.7] {$F_H$};
\addplot[thick,
    domain=0.1:0.9, 
    samples=200, 
    color=blue,
]
{2/x} node[right,pos=0.788] {$F_S$};
\addplot[only marks,mark size=0.07cm] table{
0.5 3
0.5 4
0.618034 3.23606791859
};
\end{axis}
\end{tikzpicture}
	\caption{Plots of $F_H$ in red and $F_S$ in blue}
	\label{fig:fam-h-s}
\end{figure}

Following this Talmudic judgement on the expedience of a dispute, we turn our attention to what Beit Hillel and Beit Shammai think on the matter of family planning, see Talmud, Yevamot 61b:15: 

\begin{quote}
  {\it  Mishnah: A man may not neglect the mitzva to be fruitful and multiply unless he already has children. Beit Shammai say: One fulfills this mitzva with two males, and Beit Hillel say: A male and a female, as it is stated: ''Male and female He created them'' (Genesis 5:2).}  
\end{quote}
An explanation for both requirements on how to fulfil this mitzva, a commandment to be performed as a religious duty in Judaism, can be found further in Talmud, Yevamot 61b:18:

\begin{quote}
  {\it  The Gemara asks: What is the reason of Beit Shammai? The Gemara answers: We learn this from Moses as it is written: “The sons of Moses, Gershom and Eliezer” (I Chronicles 23:15). Since Moses did not have any other children, two sons must be sufficient to fulfill the mitzva. And the reason of Beit Hillel is that we learn from the creation of the world, as mankind was created male and female. }  
\end{quote}
What are the consequences of both rabbinical interpretations for a society? More precisely, how many children would a family have on average if the society obeys the rule of Beit Hillel or Beit Shammai, respectively. This question has been addressed within a mathematical modelling by Jonathan Rosenberg in \cite{rosenberg2001}. The assumptions of his model are as follows:
\begin{itemize}
    \item[] A1: A couple begins to have children and agrees to continue until they have fulfilled their obligation from one of the aforementioned rules.
    \item[] A2: The probability of having a boy or a girl is constant each time with $p \in (0,1)$ and $1-p$, respectively,  and the genders of the children are independent events.
\end{itemize}
Let us denote the average family size according to Beit Hillel and Beit Shammai by $F_H$ and $F_S$ respectively. These quantities are not hard to represent as series, see also Section \ref{subsec:pd}:
\[
   F_H(p) = \sum_{\ell=1}^{\infty} (\ell+1)p^\ell(1-p)+\sum_{\ell=1}^{\infty} (\ell +1)p(1-p)^\ell, \quad 
   F_S(p) = \sum_{\ell=1}^{\infty} (\ell+1) p^2(1-p)^{\ell-1}.
\]
In explicit form we have:
\[
   F_H(p) = \frac{p^2 - p + 1}{p - p^2}, \quad 
   F_S(p) = \frac{2}{p}.
\]
For the equal probabilities for boys and girls, Beit Shammai outperforms Beit Hillel by exactly one child: 
\[
F_H\left(\nicefrac{1}{2}\right)=3, \quad F_S\left(\nicefrac{1}{2}\right)=4.  
\]
In this case, the asymmetric rule of Beit Shammai leads on average to a larger population than the symmetric rule of Beit Hillel. It is interesting to mention that both rules provide the same average family size if and only if the probability of having a boy corresponds to the Golden ratio, i.e.
\[
   F_H(p)=F_S(p) \quad \mbox{if and only if} \quad p = \frac{-1+\sqrt{5}}{2} \approx 0.618033\ldots
\]
For the plots of $F_H$ and $F_S$ in dependence of $p$ we refer to Figure \ref{fig:fam-h-s}.

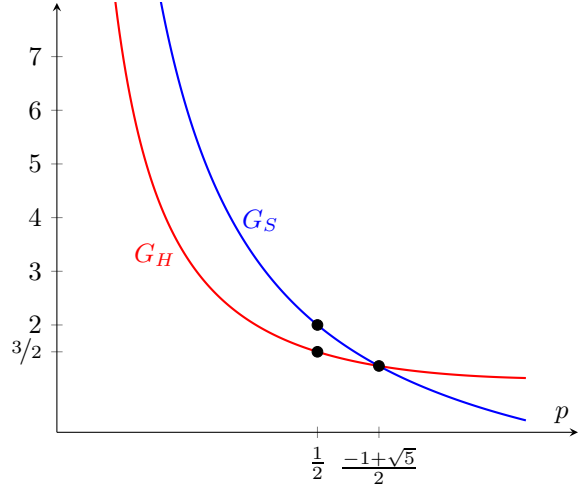
\begin{figure}[ht]
\centering
\begin{tikzpicture}[style=thick]
\begin{axis}
[
  x axis line style = {-latex},
  y axis line style = {-latex},
  axis lines = center,
  xmin=0, xmax=1, ymin=0, ymax=8,
  xlabel = {$p$},
  ylabel = {},
  xtick={0.5,0.618034}, 
  xticklabels={$\frac{1}{2}$,$\frac{-1+\sqrt{5}}{2}$},
  yticklabels={$\nicefrac{3}{2}$,2,3,4,5,6,7},
  ytick={1.5,2,3,4,5,6,7}
]
\addplot[thick,
    domain=0.1:0.9, 
    samples=200, 
    color=red,
]
{(x^2 - x + 1)/x} node[left,pos=0.7] {$G_H$};
\addplot[thick,
    domain=0.1:0.9, 
    samples=200, 
    color=blue,
]
{2*(1-x)/x} node[right,pos=0.788] {$G_S$};
\addplot[only marks,mark size=0.07cm] table{
0.5 1.5
0.5 2
0.618034 1.2360679775
};
\end{axis}
\end{tikzpicture}
	\caption{Plots of $G_H$ in red and $G_S$ in blue}
	\label{fig:girl-h-s}
\end{figure}

In this paper, we would like to go one step further compared to \cite{rosenberg2001} and to address discriminatory issues behind the family planning rules. For that, let us denote the average number of girls according to Beit Hillel and Beit Shammai by $G_H$  and $G_S$, respectively. The same applies to boys with $B_H$  and $B_S$. Analogously to above, it holds:
\[
   G_H(p) = \sum_{\ell=1}^{\infty} p^\ell(1-p)+\sum_{\ell=1}^{\infty} \ell p(1-p)^\ell, \quad 
   G_S(p) = \sum_{\ell=1}^{\infty} (\ell-1) p^2(1-p)^{\ell-1},
\]
and, explicitly:
\[
  G_H(p)= \frac{p^2 - p + 1}{p}, \quad G_S(p)=\frac{2(1 - p)}{p}.
\]
The corresponding plots are drawn in Figure \ref{fig:girl-h-s}.
Again, in the case of equal birth probabilities the average number of girls according to Beit Shammai is higher than that according to Beit Hillel:    
\[
G_H\left(\nicefrac{1}{2}\right)=\nicefrac{3}{2}, \quad G_S\left(\nicefrac{1}{2}\right)=2.  
\]
Note that this is exactly the half of the average number of children as computed before for both rules with $3$ and $4$, respectively. Hence, for the average number of boys we also obtain:
\[
  B_H\left(\nicefrac{1}{2}\right)=\nicefrac{3}{2}, \quad B_S\left(\nicefrac{1}{2}\right)=2,
\]
or in general case:
\[
  B_H(p)= F_H(p)-G_H(p)=\frac{p^2 - p + 1}{1-p}, \quad B_S(p)=F_S(p)-G_S(p)=2.
\]
What is actually striking about this observation is that not only the symmetric rule of Beit Hillel, but also the asymmetric rule of Beit Shammai discriminates on average neither girls nor boys.    
To make a precise statement, let us define the gender ratio by dividing the average number of boys through the average number of girls, respectively:
\[
   r_H(p)=\frac{B_H(p)}{G_H(p)}, \quad r_S(p)=\frac{B_S(p)}{G_S(p)}. 
\]
It is straightforward to see that the gender ratio is the same for Beit Hillel and Beit Shammai whatever the parameter $p \in (0,1)$ is taken. Namely, it always corresponds to the birth odds:
\[
   r_H(p)=r_S(p)=\frac{p}{1-p}.
\]
The goal of this paper is to show that this non-discriminatory result holds for any family planning rule requiring a certain minimum number of boys and girls.

\section{Gender ratio}
\label{gender}

We consider the family planning rules of getting children until at least $n\in \mathbb{N}\cup \{0\}$ boys and $k \in \mathbb{N}\cup \{0\}$ girls are born, cf. assumption A1. Let the births of boys and girls be independent from each other with the probability of having a boy denoted by $p \in (0, 1)$, cf. assumption A2. For short we say that a family planning $(n,k,p)$-rule is then given.  
For examle, the $(1,1,p)$-rule is due to Beit Hillel and the $(2,0,p)$-rule is due to Beit Shammai. We denote the average family size, the average number of girls, and the average number of boys by $F(n,k,p)$, $G(n,k,p)$, and $B(n,k,p)$, respectively. The gender ratio is then defined by
\[
   r(n,k,p)=\frac{B(n,k,p)}{G(n,k,p)}.
\]
In particular, 
\[
   r(1,1,p)=r_H(p), \quad r(2,0,p)=r_S(p).
\]
Now, we are ready to state the main result of the paper.

\begin{theorem}
\label{thm:main}
       The gender ratio corresponding to any family planning $(n,k,p)$-rule is equal to the birth odds, i.e. 
    \[
       r(n,k,p)=\frac{p}{1-p}.
    \]
\end{theorem}
Theorem \ref{thm:main} can be interpreted in terms of gender equality. It says that the gender ratio does not depend on a particular family planning rule implemented by the society. Since the gender ratio always corresponds to the birth odds, we conclude that the family planning rules are non-discriminatory, at least on average. The rest of this section is devoted to the alternative proofs of Theorem \ref{thm:main}.

\subsection{Mathematical induction}

As the base step, we show that for the $(0,1,p)$-rule the assertion of Theorem \ref{thm:main} holds. For this we calculate the expected number of boys:
\[
\begin{aligned}
B(0,1,p) &= \sum_{\ell=1}^\infty \ell \cdot p^{\ell} \cdot (1-p) = 
\frac p {1-p}.
\end{aligned}
\] 
Since for the expected number of girls we have $G(0,1,p)=1$, the gender ratio satisfies 
\[
r(0,1,p)= \frac {B(0,1,p)} {G(0,1,p)} = \frac {p} {1-p}.
\]
The same applies to the $(1,0,p)$-rule by a symmetry argument.

For induction step, we deduce that the $(n,k+1,p)$-rule satisfies the assertion of Theorem \ref{thm:main} by assuming that the $(n,k,p)$-rule does. Then, the symmetry of the problem also implies that the $(n+1,k,p)$-rule satisfies the assertion of Theorem \ref{thm:main} as well. 
Let $T$ denote the stopping time for the $(n,k,p)$-rule.
We distinguish two cases. First, the $(n,k+1,p)$-rule is satisfied at the time $T$, i.e. already $k+1$ girls were born up to this time. Let $p_1$ denote the probability of this case. We further denote $b_1$ and $g_1$ to be conditional expectations of the number of boys and girls, respectively, at the time $T$. Here, the expectation is conditioned on the case that the $(n,k+1,p)$-rule was already satisfied at the time $T$. Second, whenever at time $T$ only $k$ girls were born, we still have to wait until one further girl is born in order to satisfy the $(n,k+1,p)$-rule. The conditional expectation $b_3$ stands for the number of boys that are born while we wait.
For the latter we know that $b_3=B(0,1,p)$. Analogously, for the conditional expectation for girls we have $g_3=G(0,1,p)$. This is due to the independence assumed in A1.
Analogously, we denote $b_2$ and $g_2$ to be conditional expectations of the number of boys and girls, respectively, at the time $T$. However, here the expectation is conditioned on the case that  the $(n,k+1,p)$-rule is not yet satisfied at the time $T$.
Altogether, we obtain
\[
B(n,k+1,p)= p_1 b_1 + p_2 (b_2 + b_3), \qquad B(n,k,p) = p_1 b_1 + p_2 b_2,
\]
and
\[
G(n,k+1,p) = p_1 g_1 + p_2 (g_2 +g_3), \qquad G(n,k,p) = p_1 g_1 + p_2 g_2 .
\]
Since the $(0,1,p)$-rule satisfies the assertion of Theorem \ref{thm:main}, we have
\[
\frac {b_3}{g_3} = \frac {B(0,1,p)}{G(0,1,p)} = \frac p {1-p}.
\] 
By the induction hypothesis for the $(n,k,p)$-rule, it also holds:
\[
\frac {p_1 b_1 + p_2 b_2}{p_1 g_1 +p_2 g_2}    = \frac {B(n,k,p)}{G(n,k,p)} =  \frac p {1-p}.
\]
Therefore,
\[
\frac {B(n,k+1,p)}{G(n,k+1,p)} = \frac { p_1 b_1 + p_2 b_2 + p_2 b_3 }{ p_1 g_1 + p_2 g_2 +p_ 2 g_3}  = \frac {p}{1-p},
\] 
i.e. the induction step is shown.

\subsection{Doob's optional-stopping theorem}   

For the proof of Theorem \ref{thm:main} we use Doob's optional-stopping theorem from the probability theory, see e.g. \cite{william1991}.

\begin{theorem}
\label{thm:doob}
    Let $X$ be a martingale and $T$ a stopping time. If any of the conditions (i)-(iii) holds, then $X_T$ is integrable and $\mathbb{E}(X_T)=\mathbb{E}(X_0)$: 
    \begin{itemize}
        \item[] (i) $T$ is bounded;
        \item[] (ii) $X$ is bounded and $T$ is a. s. finite;
        \item[] (iii) $\mathbb{E}(T) < \infty$ and $\left| X_m - X_{m-1}\right|$ is uniformly bounded.
    \end{itemize}
\end{theorem}
In order to apply Doob's optional-stopping theorem in our setting, we construct a martingale as follows. Let is consider independent and identically distributed random variables $Y_i$, $i=1,\ldots$, given by
\[
   \mathbb{P}\left(Y_i=\frac{1}{p}\right) = p, \quad 
   \mathbb{P}\left(Y_i=-\frac{1}{1-p}\right) = 1-p.
\]
The values $\nicefrac{1}{p}$ and $-\nicefrac{1}{(1-p)}$ of $Y_i$ correspond to the birth of a boy or a girl at time $i$, respectively. Additionally, we set $Y_0=0$. Further, we define the sum of $X_i$, $i=1,\ldots, m$, as
\[
  X_m = \sum_{i=0}^{m} Y_i.
\]
Note that, by setting the random variables $B_m$ and $G_m$ as the numbers of boys or girls born up to time $m$, we obviously obtain:
\[
  X_m = \frac{B_m}{p} - \frac{G_m}{1-p}.
\]
Let now the stopping time $T$ correspond to the $(n,k,p)$-rule. By taking expectation above, we have:
\[
  \mathbb{E}(X_T) = \frac{\mathbb{E}(B_T)}{p} - \frac{\mathbb{E}(G_T)}{1-p}=\frac{B(n,k,p)}{p} - \frac{G(k,n,p)}{1-p}.
\]
Further, it is well-known in the probability theory that the sum of independent random variables with zero-mean is a martingale, see e.g. Example 10.4(a) in \cite{william1991}. Let us now check the condition (iii) from Theorem \ref{thm:doob} for the martingale $X$. For the expected stopping time to compute, consider families with $\ell+1$ children according to the $(n,k,p)$-rule, i.e. with $T=\ell +1$. Thus, each of them should consist of $n$ boys and $k$ girls. First, let us assume that a boy is born the last. Then, the births of the remaining $n-1$ boys could have happened at some of the first $\ell$ places. The number of girls here would be $\ell+1-n$. The probability of having this gender structure is $p^{n}(1-p)^{\ell+1-n}$. Second, let us assume that a girl is born the last. Then, the births of the remaining $k-1$ girls could have happened at some of the first $\ell$ places. The number of boys here would be $\ell+1-k$. The probability of having this gender structure is $p^{\ell+1-k}(1-p)^{k}$. Overall, we have:
\[
   \begin{array}{rcl}
   \mathbb{E}(T)&=& \displaystyle \sum_{\ell=n+k-1}^{\infty}
   (\ell+1) \binom{\ell}{n-1} p^{n}(1-p)^{\ell+1-n}\\ 
   &&\displaystyle + \sum_{\ell=n+k-1}^{\infty}
   (\ell+1) \binom{\ell}{k-1} p^{\ell+1-k}(1-p)^{k}.
   \end{array}
\]
Now, it is an easy task to see that the both series converge. For that, we denote
\[
    a_\ell = (\ell+1) \binom{\ell}{n-1} p^{n}(1-p)^{\ell+1-n}, \quad b_\ell = (\ell+1) \binom{\ell}{k-1} p^{\ell+1-k}(1-p)^{k}.
\]
The ratio test, see e.g. \cite{stewart2015}, applies as follows:
\[
   \lim_{\ell \rightarrow \infty} \frac{a_{\ell+1}}{a_\ell} = \lim_{\ell \rightarrow \infty} \frac{\ell+2}{\ell-n+2} (1-p) = 1-p < 1
\]
and
\[
   \lim_{\ell \rightarrow \infty} \frac{b_{\ell+1}}{b_\ell} = \lim_{\ell \rightarrow \infty} \frac{\ell+2}{\ell-k+2} p = p < 1.
\]
Altogether, it is shown that $\mathbb{E}(T) < \infty$.
Further, the differences of consecutive $X_m$'s are uniformly bounded, since
\[
  \left| X_m - X_{m-1}\right| = \left|Y_m\right| \leq \max\left\{\frac{1}{p}, \frac{1}{1-p}\right\}.
\]
Overall, Doob's optional-stopping theorem is applicable, and we get:
\[
   \mathbb{E}(X_T) = \mathbb{E}(X_0).
\]
By recalling the formula for $\mathbb{E}(X_T)$ and that $\mathbb{E}(X_0)=0$, we obtain from here:
\[
   \frac{B(n,k,p)}{p} - \frac{G(k,n,p)}{1-p} = 0.
\]
Thus, we are done and Theorem \ref{thm:main} is proved.

\subsection{Brute-force}
\label{subsec:pd}

Let us present the direct proof of Theorem \ref{thm:main}. For that, we focus on the average number of boys $B(n,k,p)$. For its computation, consider families with $\ell+1$ children according to the $(n,k,p)$-rule. Thus, each of them should consist of $n$ boys and $k$ girls. First, let us assume that a boy is born the last. Then, the births of the remaining $n-1$ boys could have happened at some of the first $\ell$ places. The number of girls here would be $\ell+1-n$. The probability of having this gender structure is $p^{n}(1-p)^{\ell+1-n}$. Second, let us assume that a girl is born the last. Then, the births of the remaining $k-1$ girls could have happened at some of the first $\ell$ places. The number of boys here would be $\ell+1-k$. The probability of having this gender structure is $p^{\ell+1-k}(1-p)^{k}$. Overall, we have:
\[
   \begin{array}{rcl}
   B(n,k,p)&=& \displaystyle \sum_{\ell=n+k-1}^{\infty}
   n \binom{\ell}{n-1} p^{n}(1-p)^{\ell+1-n}\\ 
   &&\displaystyle + \sum_{\ell=n+k-1}^{\infty}
   (\ell+1-k) \binom{\ell}{k-1} p^{\ell+1-k}(1-p)^{k}.
   \end{array}
\]
Let us simplify this expression:
\[
  \begin{array}{rcl}
   B(n,k,p)&=& \displaystyle \frac{n}{(n-1)!} p^n\sum_{\ell=n+k-1}^{\infty}
   \ell  (\ell -1) \cdots  (\ell-n+2) (1-p)^{\ell-n+1}\\ 
   &&\displaystyle + \frac{1}{(k-1)!} p(1-p)^k\sum_{\ell=n+k-1}^{\infty}
   \ell  (\ell -1) \cdots  (\ell-k+1) p^{\ell-k}.
   \end{array}
\]
We denote both series here as $S_1(n,k,p)$ and $S_2(n,k,p)$ to get overall:
\[
   B(n,k,p)= \displaystyle \frac{n}{(n-1)!} p^n S_1(n,k,p)+ \frac{1}{(k-1)!} p(1-p)^k S_2(n,k,p).
\]
Now, we derive a formula relating $S_1$ and $S_2$ to each other:
\[
  \begin{array}{rcl}
   (1-p)S_2(n,k,p)&=& \displaystyle  \sum_{\ell=n+k-1}^{\infty}
   \ell  (\ell -1) \cdots  (\ell-k+1) p^{\ell-k} 
   \\
   &&\displaystyle -\sum_{\ell=n+k-1}^{\infty}
   \ell  (\ell -1) \cdots  (\ell-k+1) p^{\ell-k+1}
   \\
   &=& \displaystyle (n+k-1)\cdots n p^{n-1} \\
   && \displaystyle +\sum_{\ell=n+k-1}^{\infty}
   (\ell+1)\ell  (\ell -1) \cdots  (\ell-k+2) p^{\ell-k+1} \\
   
   &&\displaystyle -\sum_{\ell=n+k-1}^{\infty}
   \ell  (\ell -1) \cdots  (\ell-k+1) p^{\ell-k+1} 
   \\
      &=& \displaystyle  (n+k-1)\cdots n p^{n-1} \\
     &&\displaystyle + k \sum_{\ell=n+k-1}^{\infty}
   \ell  (\ell -1) \cdots  (\ell-k+2) p^{\ell-k+1} \\
   &=& \displaystyle  (n+k-1)\cdots n p^{n-1}+ k S_1(k,n,1-p).
   \end{array}
\]
We use the latter relation to compute:
\[
\begin{array}{rcl}
   (1-p)B(n,k,p)&=& \displaystyle \frac{n}{(n-1)!} p^n(1-p) S_1(n,k,p) \\ && \displaystyle + \frac{1}{(k-1)!} p(1-p)^k (1-p)S_2(n,k,p) \\
  &=& \displaystyle \frac{n}{(n-1)!} p^n(1-p) S_1(n,k,p) \\ && \displaystyle + \frac{1}{(k-1)!} p(1-p)^k (n+k-1)\cdots n p^{n-1} \\ 
  && \displaystyle +\frac{1}{(k-1)!} p(1-p)^k k S_1(k,n,1-p) \\
  &=& \displaystyle \frac{n}{(n-1)!} p^n(1-p) S_1(n,k,p) \\ && \displaystyle + (n+k-1)\binom{n+k-2}{k-1} p^n(1-p)^k \\ 
  && \displaystyle +\frac{k}{(k-1)!} p(1-p)^k S_1(k,n,1-p).
\end{array}
\]
Due to the well-known identity
\[
  \binom{n+k-2}{k-1}=\binom{n+k-2}{n-1},
\]
the right-hand side remains invariant if substituting $n,k,p$ through $k,n,1-p$, respectively. 
Hence, the left-hand side also does, i.e.
\[
   (1-p)B(n,k,p) = p B(k,n,1-p).
\]
Obviously, the arguments of symmetry additionally provide:
\[
  G(n,k,p) = B(k,n,1-p).
\]
Overall, we thus obtained the assertion of Theorem \ref{thm:main}. The above considerations must be slightly adapted if either $n$ or $k$ is zero. Then,
$B(n,k,p)$ or $G(n,k,p)$ consists of only one of the two
series that we have computed above. Again, the similar arguments based on the hidden symmetry do the job.


\begin{remark}
Note that we can calculate the series appearing in $B(n,k,p)$ and $G(n,k,p)$ by repeatedly using the identity
$$
\begin{aligned}
    & \qquad \sum_{\ell=n+k-1}^{\infty} \ell (\ell-1) \cdots(\ell-m+1) p^{\ell-m} = \frac {d^m}{dp^m} \left( \frac {p^{n+k-1}}{1-p} \right).
\end{aligned}
$$
We thus obtain
\[
\begin{array}{rcl}
    B(n,k,p) & = & \displaystyle \frac n {(n-1)!} p^n  (-1)^{n-1} \frac {d^{n-1}}{dp^{n-1}} \left( \frac {(1-p)^{n+k-1}}{p}\right) \\ && \displaystyle + \frac {p (1-p)^k}{(k-1)!}\frac {d^{k}}{dp^{k}} \left( \frac {p^{n+k-1}}{1-p} \right).
\end{array}
\]
So the average number of boys is given by a rational function in $p$.
By recalling the identity
$
G(n,k,p) = B(k,n,1-p),
$ the gender ratio can be expressed as an explicit rational function in $p$. Nevertheless, these formulas seemed not to be rather helpful in proving that indeed 
$$
 \frac {B(n,k,p)}{G(n,k,p)} = \frac {p}{1-p}.
$$    
\end{remark}

\section{Discussion}

Let us focus here on the share of girls for the Beit Shammai rule:
\[
   g_S(p) = \frac{G_S(p)}{F_S(p)}=1-p.
\]
Note that $g_S$ relates the average number of girls to the average family size. Alternatively, we may first determine the share of girls within each family, which satisfies the Beit Shammai rule of having at least two boys, and only then we average. Let us denote this ''average share'' by $\bar g_S(p)$ and compute it in the following. For that, we consider families with $\ell+1$ children according to Beit Shammai. Thus, each of them should consist of $2$ boys and $\ell-1$ girls. Moreover, one boy is born the last, and the birth of another boy could have happened at one of the first $\ell$ places. The probability of having this gender structure is $p^2(1-p)^{\ell-1}$. Overall, we obtain
\[
  \bar g_S(p) = \sum_{\ell=1}^{\infty} \ell \frac{\ell-1}{\ell+1} p^2(1-p)^{\ell-1}. 
\]
By using the Taylor series for the logarithm, see e.g. \cite{stewart2015}, it is easy, but tedious, to compute the last expression explicitly, and we get
\[
   \bar g_S(p) = 1-2\frac{p}{1-p} \left(1+ \frac{p}{1-p} \ln p\right).
\]
In particular, for the case of equal birth probabilities we have:
\[
   \bar g_S(\nicefrac{1}{2}) = 2 \ln 2 - 1 \approx 0.386294\ldots
\]
We conclude that the average share of girls is below their birth probability:
\[
\bar g_S(\nicefrac{1}{2}) < g_S(\nicefrac{1}{2})=\nicefrac{1}{2}.
\]
In general, it can be shown that for all $p \in (0,1)$ it holds:
\[
\bar g_S(p) < g_S(p)=1-p.
\]
The latter inequality means that the average share $\bar g_S$ captures the discriminatory asymmetry in the Beit Shammai rule, whereas the share $g_S$ does not. We conclude that the Beit Shammai rule is non-discriminatory from the point of view of the society, since it focuses on $g_S$. However, from the individual point of view, where the interest is rather in $\bar g_S$, a certain discrimination is perceived.

\end{document}